\documentclass{amsart}

\usepackage{amssymb}
\usepackage[cmtip,all]{xy}
\copyrightinfo{2019}{American Mathematical Society}

\newtheorem{theorem}{Theorem}
\newtheorem{lemma}[theorem]{Lemma}
\newtheorem{corollary}[theorem]{Corollary}

\theoremstyle{definition}

\newtheorem{example}[theorem]{Example}

\theoremstyle{remark}
\newtheorem{remark}[theorem]{Remark}
\newtheorem{question}[theorem]{Question}

\numberwithin{equation}{section}

\newcommand{\VVec}{\mathbf{Vec}}
\newcommand{\Top}{\mathbf{Top}}
\newcommand{\bC}{\mathbf{C}}
\newcommand{\ZZo}{\ZZ_{\geq 0}}
\newcommand{\cI}{\mathcal{I}}
\newcommand{\del}{\partial}
\newcommand{\Add}{\mathrm{Add}}
\newcommand{\comment}[1]{}
\newcommand{\Sh}{\mathrm{Sh}}
\newcommand{\FI}{\mathbf{FI}}
\newcommand{\uV}{\underline{V}}
\newcommand{\uW}{\underline{W}}
\newcommand{\Hom}{\operatorname{Hom}}
\newcommand{\Wedge}{\bigwedge}
\newcommand{\GL}{\mathrm{GL}}
\newcommand{\ZZ}{\mathbb{Z}}
\newcommand{\NN}{\mathbb{N}}
\newcommand{\CC}{\mathbb{C}}
\newcommand{\ot}{\leftarrow}
\newcommand{\diag}{\operatorname{diag}}
\newcommand{\cha}{\operatorname{cha}}
\newcommand{\id}{\mathrm{id}}
\newcommand{\Sym}{\operatorname{Sym}}

\begin{document}

\title{Topological Noetherianity of polynomial functors}
\author{Jan Draisma}
\address{Universit\"at Bern, Mathematisches Institut, Sidlerstrasse 5,
3012 Bern; and Eindhoven University of Technology, P.O.~Box 513, 5600 MB,
Eindhoven, The Netherlands}
\email{jan.draisma@math.unibe.ch}
\thanks{The author was partially supported by the NWO Vici grant
entitled {\em Stabilisation in Algebra and Geometry}.}

\begin{abstract}
We prove that any finite-degree polynomial functor over an infinite
field is topologically
Noetherian. This theorem is motivated by the recent resolution, by
Ananyan-Hochster, of Stillman's conjecture; and a recent Noetherianity
proof by Derksen-Eggermont-Snowden for the space of cubics. Via work by
Erman-Sam-Snowden, our theorem implies Stillman's conjecture and indeed
boundedness of a wider class of invariants of ideals in polynomial rings
with a fixed number of generators of prescribed degrees.
\end{abstract}

\maketitle

\section{Introduction}

This paper is motivated by two recent developments in ``a\-symptotic
commutative algebra''. First, in \cite{Hochster16}, Hochster-Ananyan
prove a conjecture due to Stillman \cite[Problem
3.14]{Peeva09}, to the effect that the
projective dimension of an ideal in a polynomial ring generated by a fixed
number of homogeneous polynomials of prescribed degrees can be bounded
independently of the number of variables. Second, in \cite{Derksen17},
Derksen-Eggermont-Snowden prove that the inverse limit over $n$ of the
space of cubic polynomials in $n$ variables is topologically Noetherian up
to linear coordinate transformations. These two theorems show striking
similarities in content, and in \cite{Erman17}, Erman-Sam-Snowden
show that topological Noetherianity of a suitable space of tuples of
homogeneous polynomials, together with Stillman's conjecture, implies
a generalisation of Stillman's conjecture to other ideal invariants.
In addition to being similar in content, the two questions have
similar histories---e.g. both were first established for tuples of
quadrics \cite{Ananyan12, Eggermont14}---but since \cite{Hochster16}
the Noetherianity problem has been lagging behind. The goal of this
paper is to make it catch up.

\subsection{Polynomial functors} \label{ssec:Poly}

Let $K$ be an infinite field and let $\VVec$ be the category of
finite-dimensional vector spaces over $K$. We consider a covariant
polynomial functor $P:\VVec \to \VVec$ of finite degree $d$.
This means that for all $V,W$ the map $P:\Hom_\VVec(V,W) \to
\Hom_\VVec(P(V),P(W))$ is polynomial of degree at most $d$, with equality
for at least some choice of $V$ and $W$. The uniform upper bound $d$
rules out examples like $V \mapsto \Wedge^0 V \oplus \Wedge^1 V \oplus
\Wedge^2 V \oplus \ldots$.

Then $P$ splits as a direct sum $P=P_0 \oplus P_1 \oplus \cdots \oplus
P_d$ where 
\[ P_e(V):=\{p \in P(V) \mid P(t\cdot 1_V) p=t^e p \text{
for all } t\in K\}; \] 
see e.g.~\cite{Friedlander97}. For each $e$ the
map $\Hom(V,W) \to \Hom(P_e(V),P_e(W))$ is homogeneous of degree $e$,
and we have $P_d(V) \neq 0$ for all $V$ of sufficiently large dimension.

\subsection{Topological spaces over a category}

The polynomial functor $P$ will also be interpreted as a functor
to the category $\Top$ of topological spaces.  Here we equip the
finite-dimensional vector space $P(V)=P_0(V) \oplus \cdots \oplus P_d(V)$
with the Zariski-topology.
A general functor $X$ from a category $\bC$ to $\Top$ is called a
$\bC$-topological space, or $\bC$-space for short. A second $\bC$-space
$Y$ is called a $\bC$-subspace of $X$ if for each object $A$ of $\bC$
the space $Y(A)$ is a subset of $X(A)$ with the induced topology, and
if, moreover, for a morphism $f:A \to B$ the continuous map $Y(f):Y(A)
\to Y(B)$ is the restriction of $X(f)$ to $Y(A)$. We then write $Y
\subseteq X$. The $\bC$-subspace $Y$ is called closed if $Y(A)$ is
closed in $X(A)$ for each $A$; we then also call $Y$ a {\em closed
$\bC$-subset} of $X$. The $\bC$-space $X$ is called {\em Noetherian}
if every descending chain of closed $\bC$-subspaces stabilises.

Furthermore, a $\bC$-continuous map from a $\bC$-space $X$ to a
$\bC$-space $X'$ consists of a continuous map $\varphi_A:X(A) \to X'(A)$
for each object $A$, in such a way that for any morphism $f:A \to B$
we have $\varphi_B \circ X(f)=X'(f) \circ \varphi_A$. A $\bC$-homeomorphism
has the natural meaning.

\subsection{The main theorem}

We will establish the following fundamental theorem.

\begin{theorem} \label{thm:Noetherianity}
Let $\VVec$ be the category of finite-dimensional vector spaces over
an infinite field $K$, and let $P:\VVec \to \VVec$ be a finite-degree
polynomial functor. Then $P$ is Noetherian as a $\VVec$-topological space.
\end{theorem}

\begin{remark}
The restriction to infinite $K$ is crucial for the set-up and the proofs
below---e.g., it is used in the decomposition of a polynomial functors
into homogeneous parts and, implicitly, to argue that if an algebraic
group $\GL_n(K)$ preserves an ideal, then so does its Lie algebra. In
future work, we will pursue versions of Theorem~\ref{thm:Noetherianity}
over $\ZZ$ and possibly over finite fields. 
\end{remark}

Theorem~\ref{thm:Noetherianity} will be useful in different contexts where
finiteness results are sought for. In the remainder of this section we
discuss several such consequences; since an earlier version of this
paper appeared on the arxiv, several other ramifications
have appeared, e.g., in \cite{Erman17}. 

\subsection{Equivariant Noetherinity of limits} The polynomial functor
$P$ gives rise to a topological space $P_\infty:=\lim_{\ot n} P(K^n)$,
the projective limit along the linear maps $P(\pi_n): P(K^{n+1}) \to
P(K^n)$ where $\pi_n$ is the projection $K^{n+1} \to K^n$ forgetting
the last coordinate. By functoriality, each $P(K^n)$ is acted upon by 
the general linear group $\GL_n$, the map $P(\pi_n)$ is
$\GL_n$-equivariant if we embed $\GL_n$ into $\GL_{n+1}$ via $g \mapsto
\diag(g,1)$, and hence $P_\infty$ is acted upon by the
direct limit $\GL_\infty:=\bigcup_n \GL_n$, the group of all invertible
$\NN \times \NN$-matrices which in all but finitely many entries equal
the infinite identity matrix.

Given a closed $\VVec$-subset
$X$ of $P$, the inverse limit $X_\infty:=\lim_{\ot n} X(K^n)$ is a closed,
$\GL_\infty$-stable subset of $P_\infty$, and using embeddings $K^n \to
K^{n+1}$ appending a zero coordinate, one finds that $X_\infty$ surjects
onto each $X(K^n)$. Conversely, given a closed, $\GL_\infty$-stable subset
$Y$ of $P_\infty$, then for any finite-dimensional vector space $V$ and
any linear isomorphism $\varphi:K^n \to V$, we set $X(V):=P(\varphi)(Y_n)$,
where $Y_n$ is the image of $Y$ in $P(K^n)$. 

One can check that $V \mapsto
X(V)$ is a closed $\VVec$-subset of $P$ and that this
construction $Y \mapsto X$ is
inverse to the construction $X \mapsto X_\infty$ above. Thus the theorem
is equivalent to the following corollary. 

\begin{corollary}
Let $\VVec$ be the category of finite-dimensional vector spaces over an
infinite field $K$, let $P:\VVec \to \VVec$ be a finite-degree polynomial
functor, and equip $P_\infty:=\lim_{\ot n} P(K^n)$ with the inverse-limit
topology of the Zariski topologies on the $P(K^n)$. Then $P_\infty$ is
{\em $\GL_\infty$-Noetherian}, i.e., every chain 
$P_\infty \supseteq Y_1 \supseteq Y_2 \supseteq \ldots$ 
of $\GL_\infty$-stable closed subsets of $P_\infty$ is eventually
constant. Equivalently, every $\GL_\infty$-stable closed
subset $Y$ of $P_\infty$ is the set of common zeroes of finitely many
$\GL_\infty$-orbits of polynomial equations. 
\end{corollary}

\begin{example}
The paper \cite{Derksen17} concerns the case where $P=S^3$, the third
symmetric power. In this case, $P_\infty$ is the space of infinite cubics
$\sum_{1 \leq i \leq j \leq k} c_{ijk} x_i x_j x_k$, and $\GL_\infty$
acts by linear transformations that affect only finitely many of the
variables $x_i$. \hfill $\clubsuit$
\end{example}

\begin{remark}
The proofs below could have been formulated directly in this
infinite-dimensional setting, rather than the finite-dimensional,
functorial setting. However, one of the key techniques, namely, shifting
a polynomial functor by a constant vector space, is best expressed in
the functorial language. Moreover, the functorial language allows us to
stay in the more familiar realm of finite-dimensional algebraic geometry.
\end{remark}

\subsection{Generalisations of Stillman's conjecture}
In \cite{Erman17}, Erman, Sam and Snowden use the following 
special case of Theorem~\ref{thm:Noetherianity}.

\begin{corollary} \label{cor:S}
Let $K$ be an infinite field, fix natural numbers $d_1,\ldots,d_k$, and
consider the polynomial functor $P:V \mapsto S^{d_1} V \oplus \cdots
\oplus S^{d_k} V$. Then $P$ is a Noetherian $\VVec$-topological space,
and hence its limit $P_\infty$ is $\GL_\infty$-Noetherian.
\end{corollary}

Let $\mu$ be a function that associates a number $\mu(I) \in \ZZ \cup
\{\infty\}$ to any homogeneous ideal $I$ in a symmetric algebra $S V$ on
$V \in \VVec$, in such a way that $\mu(\langle S(\varphi)I \rangle)=\mu(I)$
for any injective linear map $\varphi:V \to W$ with induced homomorphism
$S(\varphi):S V \to S W$, and such that $\mu$ is upper semicontinuous in
flat families. In \cite{Erman17} the following is proved.

\begin{theorem}[\cite{Erman17}]
Corollary~\ref{cor:S} implies that for any ideal invariant $\mu$
with the properties above there exists a number $N$ such that for all
$V \in \VVec$, any ideal $I \subseteq S V$ generated by $k$ homogeneous
polynomials of degrees $d_1,\ldots,d_k$
either has $\mu(I) \leq N$ or $\mu(I)=\infty$.
\end{theorem}

The crucial point here is that $N$ does not depend on $\dim V$. Stillman's
conjecture \cite[Problem 3.14]{Peeva09} is this statement for $\mu$
equal to the projective dimension, and it is used in the proof of the
generalisation just stated. However, in a follow-up paper \cite{Erman18},
the same authors give two new proofs of Stillman's conjecture, one of
which uses Corollary~\ref{cor:S}.  An algorithmic variant of this latter
proof is presented in \cite{Draisma18a}.

\subsection{Twisted commutative algebras} 
For $K=\CC$, the algebra of polynomial
functions on $P_\infty$, i.e., the direct limit of the symmetric algebras
on $P_n(K^n)^*$, is a {\em twisted commutative algebra} in one of its
incarnations \cite{Sam12,Sam12b}. In this context,
Theorem~\ref{thm:Noetherianity} says the following.

\begin{corollary}
Any finitely generated twisted commutative algebra over $\CC$ 
is topologically Noetherian.
\end{corollary}

\subsection{Functors from $\VVec^\ell$} Theorem~\ref{thm:Noetherianity}
has the following generalisation to functors that take several distinct
vector spaces as input.

\begin{corollary} \label{cor:Vecl}
Let $K$ be an infinite field, $\VVec$ the category of finite-dimensional
vector spaces over $K$, $\ell$ a positive integer, and $P$ a functor
from $\VVec^\ell$ to $\VVec$ such that for any $\uV, \uW \in \VVec^\ell$ the map
$\Hom_{\VVec^\ell}(\uV,\uW) \to \Hom_\VVec(P(\uV),P(\uW))$ is polynomial of
uniformly bounded degree. Then $P$ is a Noetherian $\VVec^\ell$-topological
space. 
\end{corollary}

Note that the group of automorphisms of $(V_1,\ldots,V_\ell)$ is $\prod_i
\GL(V_i)$, which when the $V_i$ are all equal contains a diagonal copy
of $\GL(V)$. This suggests that Theorem~\ref{thm:Noetherianity} for $V
\mapsto P(V,\ldots,V)$ is in fact stronger than this corollary, as we
prove now.

\begin{proof}[Proof of Corollary~\ref{cor:Vecl} from
Theorem~\ref{thm:Noetherianity}.]
Let $X_1 \supseteq X_2 \supseteq \cdots$ be a chain of closed
$\VVec^\ell$-subsets of $P$. Let $Q:\VVec \to \VVec$ be the functor that sends $V$ to
$P(V,\ldots,V)$, and set $Y_n(V):=X_n(V,\ldots,V)$, a closed $\VVec$-subset
of the polynomial functor $Q$. By Theorem~\ref{thm:Noetherianity}, the
sequence $(Y_n)_n$ is constant for $n$ at least some $n_0$. We claim that
so is $(X_n)_n$. Indeed, let $\uV=(V_1,\ldots,V_\ell) \in \VVec^\ell$.
Choose a $V \in \VVec$ with $\dim V \geq \dim V_i$ for each $i$, and choose
surjections $\pi_i:V \to V_i$ and injections $\iota_i:V_i \to V$ such
that $\pi_i \circ \iota_i=1_{V_i}$. Then for $n \geq n_0$ we have
\begin{align*}
X_n(V_1,\ldots,V_\ell) &= P(\pi_1,\ldots,\pi_\ell)
P(\iota_1,\ldots,\iota_\ell) X_n(V_1,\ldots,V_\ell) \\
&\subseteq P(\pi_1,\ldots,\pi_\ell) X_n(V,\ldots,V) \\
&= P(\pi_1,\ldots,\pi_\ell) X_{n+1}(V,\ldots,V) 
\subseteq X_{n+1}(V_1,\ldots,V_\ell),
\end{align*}
as desired.
\end{proof}

\subsection{Slice rank} Taking $P(V_1,\ldots,V_\ell)=V_1 \otimes \cdots \otimes
V_\ell$, a tensor in $P(\uV)$ is said to have slice rank $1$ if it is
nonzero and of the form $v \otimes A$ for a $v$ in one of the $V_i$ and
an $A \in \bigotimes_{j \neq i} V_j$. A tensor has slice rank at most $k$
if it is a sum of at most $k$ tensors of slice rank $1$; in this sum the
the slice index $i$ may vary through $\{1,\ldots,d\}$.  Being of slice
rank at most a fixed number $k$ is a Zariski-closed condition (see Tao
and Sawin's blog post \cite{Tao16}) and preserved under tensor products
of linear maps.  Corollary~\ref{cor:Vecl} implies the following.

\begin{corollary}
Let $\VVec$ be the category of finite-dimensional vector spaces over an
infinite field $K$. For fixed $\ell$ and $k$, there exists a tuple $\uW
\in \VVec^\ell$ such that for all $\uV \in \VVec^\ell$ a tensor $T \in V_1
\otimes \cdots \otimes V_\ell$ has slice rank at most $k$ if and only if
for all $\ell$-tuples of linear maps $\varphi_i:V_i \to W_i$ the tensor
$(\varphi_1 \otimes \cdots \otimes \varphi_\ell)T$ has slice rank at most $k$.
\end{corollary}

Equivalently, in the space of infinite by infinite by \ldots by infinite
$\ell$-way tensors, the variety of slice-rank at most $k$ tensors is
defined by finitely many $\GL_\infty^\ell$-orbits (and even finitely many
$\GL_\infty$-orbits, acting diagonally) of polynomial equations.  A more
in-depth study of the algebraic geometry of slice rank is forthcoming
work with Oosterhof.

\subsection{Related work}

Theorem~\ref{thm:Noetherianity} fits in a trend at the interface between
representation theory, algebraic geometry, commutative algebra, and
applications, which studies algebraic structures over some base category
and aims to establish stabilisation results. Recent examples, in addition
to those referenced above, include the theory of modules over the category
$\FI$ of finite sets with injective maps \cite{Church12}; Gr\"obner
techniques \cite{Sam14} for modules over more general combinatorial
categories that, among other things, led to a resolution of the
Artinian conjecture \cite{Putman14} and to a resolution of a conjecture
by Rauh-Sullivant \cite{Rauh14} on iterated toric fibre products
\cite{Draisma16b}; and finiteness results for secant varieties of Segre
and Segre-Veronese embeddings; see \cite{Landsberg04,Raicu10,Sam15b,Sam16}
and the notion of inheritance in \cite{Landsberg11}. The current paper,
while logically independent of these results, was very much influenced
by the categorical viewpoint developed in these papers.

\subsection*{Acknowledgments}

I thank Arthur Bik, Micha\l\ Laso\'n, Florian Oosterhof, and Andrew
Snowden for useful discussions and comments on an earlier version
of this paper. I also thank the organisers of the April 2016 Banff
workshop on {\em Free Resolutions, Representations, and Asymptotic
Algebra} for bringing together participants with a wide variety of
backgrounds---they have strongly influenced my understanding of these
infinite algebraic structures. Finally, I thank the referees for several
valuable suggestions, including the running example~\ref{ex:Running}.

\section{Proof of the main theorem} 

\subsection{Overview of the proof}
The proof of Theorem~\ref{thm:Noetherianity} is a double induction. The
{\em outer induction} is on the polynomial functor $P$ via a
(``lexicographic'') partial order $\prec$ on the class of polynomial
functors introduced in \S \ref{ssec:Order}. Using classical work by
Friedlander-Suslin, we prove that this is a well-founded order. Any 
degree-zero polynomial functor, i.e., a constant vector space $U$
independent of the input $V$, is smaller than all polynomial functors
of positive degree, and Hilbert's basis theorem yields the base case of
the induction.

So when we want to prove the theorem for $P$, we may assume that it
holds for all polynomial functors smaller than $P$. We then show that
every closed $\VVec$-subset $X$ of $P$ is Noetherian, by the {\em inner
induction} on the smallest degree of a nonzero equation $f \in K[P(U)]$
vanishing on $X(U)$ for some $U \in \VVec$. Roughly, this works as follows
(see the next paragraph for subtleties): fix an irreducible component $R$
of the highest-degree part $P_d$ of $P$, and find a direction $r_0 \in
R(U)$ such that the directional derivative
$h:={\partial{f}}/{\partial
r_0}$ is not identically zero. Let $Y$ be the largest $\VVec$-closed
subset of $P$ on which $h$ vanishes identically. Since $h$ has lower
degree than $f$, $Y$ is Noetherian by the inner induction hypothesis. On
the other hand, set $P'(V):=P(U \oplus V)$ and $Q'(V):=P'(V)/R(V)$, so
that $Q' \prec P$; this is discussed in \S\ref{ssec:Shift}. 
In \S\ref{ssec:Proof} we show that the complement $Z$ of $Y$ has a closed
embedding into a basic open subset of $Q'$, so $Z$ 
Noetherian by the outer induction hypothesis. Hence $X$, the union of
two Noetherian spaces, is Noetherian.

There are four subtleties: First, $f$ may not depend on the coordinates
on $R(U)$, so that ${\partial{f}}/{\partial r_0}=0$ for all $r_0 \in
R(U)$. We therefore need to look for $f$ in the ideal of $X$ that are
nonzero even after modding out the ideal of the projection of $X$ in
$Q:=P/R$; see \S\ref{ssec:Splitting}. Second, in positive characteristic,
directional derivatives (linearisations) do not necessarily behave well;
we replace these by additive polynomials in \S\ref{ssec:Additive}. Third,
the closed embedding is in fact a Zariski homeomorphism to a closed
subset; see \S\ref{ssec:Closed} for the relevant lemma.  Fourth, it is not
quite $Z$ that embeds into a basic open subset of $Q'$---$Z$
is not functorial in $V$---but rather the locus $Z'(V)$ in $X'(V):=X(U
\oplus V)$ where $h$ (which involves only coordinates on the constant
vector space $U$) is nonzero. The closed embedding is then just the
restriction of the projection $P'(V) \to Q'(V)$ along $R(V)$. Smearing
around $Z'(V)$ by $\GL(U \oplus V)$, we obtain $Z(U \oplus V)$, and in
\S\ref{ssec:Proof} we show that this is good enough.

\begin{example} \label{ex:Running}
As a running example to illustrate the proof, we assume that $\cha
K \neq 2$, $P(V)=V \otimes V = Q(V) \oplus R(V)$ where $Q(V)$ is the space of
symmetric two-tensors (matrices) and $R(V)$ is the
space of skew-symmetric tensors. Let
$(x_{ij})_{ij}$ be the standard coordinates on $P(K^n)$,
$(y_{ij})_{i \leq j}$
be the standard coordinates on $Q(K^n)$ extended to $P(K^n)$ by declaring
them zero on $R(K^n)$, and let $(z_{ij})_{i<j}$ be the standard
coordinates on $R(K^n)$, similarly extended to $P(K^n)$. 

We then have
\[ x_{ij}=\begin{cases}
y_{ij}+z_{ij} & \text{if }i<j;\\
y_{ji}-z_{ji} & \text{if $j<i$; and} \\
y_{ii} & \text{if }i=j. 
\end{cases}
\]
Take $X(V)=\{v \otimes w \mid v,w \in V\}$, the variety of
rank-one tensors. Then $X(K^1)$ is the entire space $P(K^1)$,
but $X(K^2) \subsetneq P(K^2)$, so we may take $U=K^2$ and
for $f$ the determinant 
\[
f=x_{11}x_{22}-x_{12}x_{21}=y_{11}y_{22}-(y_{12}+z_{12})(y_{12}-z_{12})
=y_{11}y_{22}-y_{12}^2+z_{12}^2.
\]
We take $r_0:=e_1 \otimes e_2 - e_2 \otimes e_1$ and find
\[ 
h=\partial f / \partial r_0=\partial f / \partial z_{12} = 2 z_{12}.
\]
In this case, for $n \geq 2$, $Y(K^n)$ is the subvariety of $X(K^n)$ on
which the $\GL_n$-orbit of $z_{12}$ vanishes identically, i.e., $Y(V)$
is the set of rank-one tensors in $Q(V)$. This is a coincidence; in the general
setting of the proof, $Y(V)$ does not embed into $Q(V)$, but it always has a
lower-degree polynomial vanishing on it. We
discuss the complement $Z(V)$ in Example~\ref{ex:Running4}. \hfill
$\clubsuit$
\end{example}

\subsection{A well-founded order on polynomial functors}
\label{ssec:Order}

We will prove Theorem~\ref{thm:Noetherianity} by induction on
the polynomial functor, along a partial order that we introduce now.
Define a relation $\prec$ on polynomial functors of finite degree by $Q
\prec P$ if $Q \not \cong P$ and moreover if $e$ is the highest degree
with $Q_e \not \cong P_e$, then $Q_e$ is a homomorphic image of $P_e$;
this is a partial order on (isomorphism classes of)
polynomial functors.

\begin{lemma} \label{lm:WellFounded}
The relation $\prec$ is a well-founded order on 
polynomial functors of finite degree.
\end{lemma}

\begin{proof}
It suffices to prove that this order is well-founded when restricted to
polynomial functors of degree at most a fixed number $d$. By
\cite[Lemma 3.4]{Friedlander97},
if $V$ is any vector space of dimension at least $d$, then the map $P \mapsto
P(V)$ is an equivalence of abelian categories from polynomial functors of degree
at most $d$ and finite-dimensional polynomial $\GL(V)$-representations of 
degree at most $d$. Hence $Q
\prec P$ implies that the sequence $(\dim Q_e(V))_{e=1}^d$ is strictly smaller
than the sequence $(\dim P_e(V))_{e=1}^d$ in the lexicographic order
where position $e$ is more significant than position $e-1$. Since this
lexicographic order is a well-order, $\prec$ is well-founded.
\end{proof}

\subsection{$\VVec$-varieties and their ideals}
\label{ssec:Var}

Write $K[P]$ for the contravariant functor from $\VVec$ to $K$-algebras
that assigns to $V$ the coordinate ring $K[P(V)]$.  A closed $\VVec$-subset
$X$ of $P$ will be called a $\VVec$-variety in $P$, and denoted $X
\subseteq P$. Its ideal is a contravariant functor that sends $V$ to
the ideal of $X(V)$ inside $K[P(V)]$.

Using scalar multiples of the identity $V \to V$ and the fact that $K$
is infinite, one finds that the ideal of $X(V)$ is homogeneous with respect
to the $\ZZo$-grading that assigns to the coordinates on $K[P_e(V)]$
the degree $e$. The degree function $\deg$ on $K[P(V)]$
and its quotients by homogeneous ideals is defined
relative to this grading.

\begin{example} \label{ex:Running2}
In our running Example~\ref{ex:Running}, $f,h$ have degrees $4,2$,
respectively. \hfill $\clubsuit$
\end{example}

\subsection{The shift functor}
\label{ssec:Shift}

Fixing a $U \in \VVec$, we let $\Sh_U:\VVec \to \VVec$ be the {\em shift
functor} that sends $V \mapsto U \oplus V$ and the homomorphism $\varphi:V
\to W$ to the homomorphism $\Sh_U(\varphi):=\id_U \oplus \varphi:
\Sh_U(V)
\to \Sh_U(W)$. 

\begin{lemma} \label{lm:Shift}
For any polynomial functor $P$ of degree $d$, $P \circ \Sh_U$
is a polynomial functor of degree $d$ whose degree-$d$ 
homogeneous part is canonically isomorphic to that of $P$.
\end{lemma}

\begin{proof}
Set $P':=P \circ \Sh_U$. For $V,W \in \VVec$ the map $\Hom(V,W) \to
\Hom(P'(V),P'(W))$ given by $\varphi \mapsto P'(\varphi)$ is the composition of
the affine-linear map $\varphi \mapsto 1_U \oplus \varphi$ and the polynomial
map $\psi \mapsto P(\psi)$ of degree at most $d$, so $P'$ is a polynomial
functor of degree at most $d$.

For $V \in \VVec$ let $\iota_V$ be the embedding $V \to U \oplus V,\
v \mapsto (0,v)$ and let $\pi_V$ be the projection $U \oplus V \to V,\ (u,v)
\mapsto v$. These give rise to morphisms of polynomial functors $\alpha:P
\to P'$ and $\beta:P' \to P$ given by $\alpha(V):=P(\iota_V):P(V) \to
P(U \oplus V)$ and $\beta(V):=P(\pi_V)$. Straightforward computations
show that $\alpha$ and $\beta$ map each homogeneous part $P_e$ into
$P'_e$ and vice versa, and that $\beta \circ \alpha$ is the identity.
Conversely, for $q$ in the highest-degree part $P_d'(V)$ we have
$P(1_U \oplus t 1_V)q=t^d q$ for all $t \in K$. The
coefficient of $t^d$ in the left-hand side
equals $P(0 \oplus 1_V)q$, so we have $P(0
\oplus 1_V)q=q$ and therefore 
\[ \alpha(V) \beta(V) q = P(\iota_V) P(\pi_V) q
= P(0 \oplus 1_V) q = q, \]
which proves that $\beta:P'_d(V) \to P_d(V)$ is indeed a
linear isomorphism.
\end{proof}

\begin{example}
If $P(V)=S^d V$, then $(P \circ \Sh_U)(V)=S^d (U \oplus V) \cong
\bigoplus_{e=0}^d S^{d-e} U \otimes S^e V$, so $P \circ \Sh_U$ equals $P$
plus a polynomial functor of degree $d-1$. \hfill
$\clubsuit$
\end{example}

\begin{example} \label{ex:Running3}
In our running Example~\ref{ex:Running}, $P(V)=V \otimes
V=Q(V) \oplus R(V)$, $U=K^2$ and 
\begin{align*} P \circ \Sh_U(V)&=(U \oplus V) \otimes (U \oplus V)=U \otimes U
+ U \otimes V + V \otimes U + V \otimes V\\ &\cong K^4 +
V^4 + Q(V) + R(V); 
\end{align*}
note that the degree-$2$ part of $(P \circ
\Sh_U)/R$ is $Q$, so that $(P \circ \Sh_U)/R \prec P$.
\hfill $\clubsuit$
\end{example}

\subsection{Splitting off a term of highest degree}
\label{ssec:Splitting}

Assume that $P$ is a polynomial functor of degree $d>0$.  Let $R$
be any irreducible sub-polynomial functor of the highest-degree part
$P_d$ of $P$, define $Q:=P/R$, and let $\pi:P \to Q$ be the natural
projection. Then $K[Q]$ embeds into $K[P]$ via the pull-back of $\pi$. If
$X$ is a $\VVec$-variety in $P$, then let $X_Q$ be the $\VVec$-variety in
$Q$ defined by setting $X_Q(V)$ equal to the Zariski-closure in $Q(V)$
of $\pi(V)(X(V))$.

We will think of $X$ as a $\VVec$-variety over $X_Q$. Accordingly, we
write $\cI_X$ for the contravariant functor that assigns to $V$ the
ideal of $X(V)$ in $K[\pi(V)^{-1}(X_Q(V))]$, the quotient of $K[P(V)]$
by the ideal in $K[P(V)]$ generated by the ideal of
$X_Q(V)$ in $K[Q(V)]$. 

In particular, we have $\cI_X=0$ if and only if for all $V$ we have
$X(V)=\pi^{-1}(X_Q(V))$. We write $\delta_X \in \{1,2,\ldots,\infty\}$ for
the minimal degree of a nonzero homogeneous polynomial $f \in \cI_X(V)$ as
$V$ runs over $\VVec$. Note that any polynomial in $K[P(V)]$ of degree $0$
is contained in $K[P_0(V)] \subseteq K[Q(V)]$; here we use that $d>0$. So
if a degree-$0$ polynomial vanishes on $X(V)$, then it is an equation for $X_Q(V)$ and
has already been modded out in the definition of $\cI_X$. This explains
why $\delta_X \geq 1$.  Furthermore, note that $\delta_X=\infty$ if and
only if $\cI_X=0$.

\begin{example}
In Examples~\ref{ex:Running},\ref{ex:Running2}, $\delta_X=\deg f=4$.
\hfill $\clubsuit$
\end{example}

\subsection{Additive polynomials as directional derivatives}
\label{ssec:Additive}

For a finite-dimensional vector space $W$ over the infinite field $K$,
we write $\Add(W)$ for the subset of $K[W]$ consisting of polynomials $f$
such that $f(v+w)=f(v)+f(w)$ for all $v,w \in W$. Then $\Add(W)$ is a
$K$-subspace of $K[W]$, and equal to $W^*$ when $\cha K=0$. In general, if
we let $p$ be the characteristic exponent of $K$---so $p=1$ if $\cha K=0$
and $p=\cha K$ otherwise---and if we choose a basis $x_1,\ldots,x_n$ of
$W^*$, then $\Add(W)$ has as a basis the polynomials $x_i^{p^e}$ where $i$
runs through $\{1,\ldots,n\}$ and $e$ through $\ZZo$ if $p>1$ and through
$\{0\}$ if $p=1$. The span of these for fixed $e$ is denoted $\Add(W)_e$.

\begin{lemma} \label{lm:Hasse}
Let $W' \supseteq W$ be finite-dimensional vector spaces over the
infinite field $K$, let $f$ be a polynomial on $W'$, and consider the
expression $f(w'+tw)$, a polynomial function of the triple $(w',t,w)
\in W' \times K \times W$. Then one of the following hold:
\begin{enumerate}
\item $f(w'+tw)$ is independent of $t$; this happens if and
only if $f$ factors through the projection $\pi_{W'/W}:W' \to W'/W$; or
\item the nonzero part of $f(w'+tw)$ of lowest degree in $t$ is of the form 
$t^{p^{e}} h(w',w)$ for a unique $e$
(taken $0$ if $\cha K=0$). Then for
each fixed $w' \in W'$ the map $w \mapsto h(w',w)$ is in $\Add(W)_e$.
\end{enumerate}
\end{lemma}

\begin{proof}
For $w \in W$, let $D_w^{(r)}:K[W'] \to K[W']$ denote the $r$-th
Hasse directional derivative in the direction $w$. This
linear map is defined in
characteristic $0$ by $D_w^{(r)} g =\frac{1}{r!}(\frac{\partial}{\partial
w})^r g$ and in arbitrary characteristic by realising that the latter
expression actually has integer coefficients relative to any monomial
basis of $K[W']$. Explicitly: let $x_1,\ldots,x_{n-1}$ be a basis of
$(W'/K w)$ and let $x_n \in (W')^*$ with $x_n(w)=1$. Then
\[ D_w^{(r)} x_1^{a_1} \cdots x_n^{a_n}:=\binom{a_n}{r}
x_1^{a_1} \cdots x_{n-1}^{a_{n-1}} x_n^{a_n-r}; \]
in particular, if $r=a_n$, then this is nonzero as a
polynomial over $K$, even when $a_n \cdots (a_n-r+1)$ is zero. A straightforward
check shows that this is independent of the choice of basis, and that
$D_{cw}^{(r)}=c^r D_w^{(r)}$.

Taylor's formula in arbitrary characteristic reads
\[ f(w'+tw)=\sum_{r \geq 0} (D_w^{(r)} f)(w') t^r. \]
If $D_w^{(r)} f=0$ for all $r>0$, then we see from the above that $f$
does not involve the variable $x_n$, i.e., $f$ factors through $W'
/\langle w \rangle$. Similarly, if $D_w^{(r)}f=0$ for all $w \in W$
and all $r>0$, then then $f$ factors through $W'/W$.

Suppose that there exist $r>0$ and $w \in W$ such that $D_w^{(r)}f \neq
0$; take such a pair $(r,w)$ with $r$ minimal. Then, in the coordinates
above, $f$ contains a monomial $x_1^{a_1} \cdots x_n^{a_n}$ for
which $\binom{a_n}{r}$ is nonzero in $K$, but
$\binom{a_n}{r'}=0$ in $K$ 
for all $r'$ with $0<r'<r$. By Lucas's theorem on binomial
coefficients, $r$ is a power of $p$, $a_n$ is divisible by $r$, and
$\binom{a_n}{r}=a_n/r$ in $K$. Since in fact $\binom{a_n'}{r'}=0$
holds for all $r'<r$ and {\em all} monomials $x_1^{a_1'} \cdots
x_n^{a_n'}$ in $f$, $a_n'$ is a multiple of $r$, and hence
$f=g(x_1,\ldots,x_{n-1},x_n^r)$ for a unique polynomial
$g(x_1,\ldots,x_{n-1},y_n)$. Moreover, $D_w^{(r)} f=(\partial g/\partial
y_n)(x_1,\ldots,x_{n-1},x_n^r)$.

More generally, let $x_1,\ldots,x_k$ be a basis of $(W'/W)^*$
and extend to a basis $x_1,\ldots,x_n$ of $(W')^*$. 

We then find, by minimality of $r$, that 
$f=g(x_1,\ldots,x_k,x_{k+1}^r,\ldots,x_n^r)$ for a unique polynomial
$g(x_1,\ldots,x_n,y_{k+1},\ldots,y_n)$, and for all $u \in
W$ we have 
\[ D_{u}^{(r)} f=\sum_{i=k+1}^n (x_i(u))^r (\partial g/\partial
y_i)(x_1,\ldots,x_k,x_{k+1}^r,\ldots,x_n^r). \]
Since $r=p^e$ for some $e$, the right-hand side is additive
in $u$,
which concludes the proof of the lemma.
\end{proof}

For $W' \supseteq W$ and $f \in K[W']$ and $h \in K[W'] \otimes K[W]$ as
in the second case of the lemma and for $w \in W$ we write $\del_w f \in
K[W']$ for the polynomial $w' \mapsto h(w',w)$, and call this the {\em
directional derivative} of $f$ in the direction $w$.  This polynomial
has degree less than $\deg f$, and agrees with the usual directional
derivative for $\cha K=0$. Note that $\del_w f$ depends on the choice
of $W$ inside $W'$, not just on $w$: if $f(w'+tw)$ depends on $t$ but
vanishes to a higher degree at $t=0$ for a specific $w \in W$ than it does
for general $w \in W$, then we have $\del_w f=0$. 

Also in the first case
of the lemma we write $\del_w f:=0$ for all $w \in W$. We extend the
notation to rational functions with nonzero denominator $h \in K[W'/W]$
by $\del_w(f/h):=(\del_w f)/h$. The following lemma is immediate from
Lemma~\ref{lm:Hasse}.

\begin{lemma}
For $W' \supseteq W$, $f \in K[W']$, $e \in \ZZ_{\geq 0}$ as
in Lemma~\ref{lm:Hasse}, $h \in K[W'/W] \setminus \{0\}$, 
we have $\del_{v+w} (f/h)= \del_v
(f/h) + \del_w (f/h)$
and $\del_{cw} (f/h)=c^{p^e} \del_w (f/h)$ for $c \in K$ and $w \in W$. \hfill
$\square$
\end{lemma}

\begin{example}
Assume that $p>2$. Let $W'=K^3$ with standard coordinates $x,y,z$, let $W$
be the span of the first two standard basis vectors, and let $f=y^{p^2}
z^2 + x^{2p} y^{p^2} z$. Then
\begin{align*} f((x,y,z)+t(a,b,0))&=(y+tb)^{p^2} z^2 + (x+ta)^{2p}
(y+tb)^{p^2} z \\
&= f(x,y,z) + t^p (2 a^p x^p y^{p^2} z) + \cdots \end{align*}
where the remaining terms are divisible by $t^{2p}$. Hence
$\del_{(a,b,0)} f = 2 a^p x^p y^{p^2} z$. 
\hfill
$\clubsuit$
\end{example}

\subsection{A closed embedding}
\label{ssec:Closed}

Retaining the notation in the previous section, let $B$ be a basic open
subset in $W'/W$ defined by the non-vanishing of some polynomial $h
\in K[W'/W] \setminus \{0\}$ (we allow $h=1$, in which case $B=W'/W$).
Let $Z$ be a Zariski-closed subset of $A:=\pi_{W'/W}^{-1}(B) \subseteq
W'$ and let $J$ be the ideal of $Z$ inside $K[A]=K[W'][1/h]$. Fix a number $e \in
\ZZ_{\geq 0}$, equal to $0$ if $\cha K=0$, and let $J_e$ be the set of
elements $k \in J$ such that $k(a+tw)=k(a)+t^{p^e}(\del_w k)(a)$ for all
$a \in A, t \in K, w \in W$ (so $k$ is {\em affine-additive} in $W$ with
additive part of degree $p^e$). Note that via the pull-back $\pi_{W'/W}^*:
K[B] \to K[A]$, $J_e$ is a $K[B]$-submodule of $K[A]$.

\begin{lemma} \label{lm:Closed}
Assume that $K$ is algebraically closed and suppose that for each $a \in
A$ the map $J_e \to \Add(W)_e, k \mapsto (w \mapsto (\partial_w k)(a))$
is surjective. Then $\pi_{W'/W}$ restricts to a Zariski-homeo\-morphism
from $Z$ to a closed subset of $B$.
\end{lemma}

\begin{proof}
Fix any tuple $x_1,\ldots,x_n \in (W')^*$ whose restrictions to $W$ form
a basis of $W^*$. Then the natural map $K[B][x_1,\ldots,x_n] \to K[A]$
is an isomorphism by which we identify the two algebras. Under this
identification, each element of $J_e$ can be written uniquely as $k_0 +
\sum_{j=1}^n k_j x_j^{p^e}$ for suitable elements $k_0,k_1,\ldots,k_n
\in K[B]$. Let $k^{(1)},\ldots,k^{(m)}$ be $K[B]$-module generators of
$J_e$, and let $k^{(i)}_j \in K[B]$ be the coefficient of $x_j^{p^e}$
in $k^{(i)}$. Let $M \in K[B]^{m \times n}$ be the matrix whose $(i,j)$
entry equals $k^{(i)}_j$. The condition in the lemma says that $M(b)$
has rank $n$ for all $b \in B$.

Since $K$ is algebraically closed,
the Nullstellensatz implies that $1$ lies in the ideal generated by the
nonzero $n \times n$-subdeterminants $\Delta_1,\ldots,\Delta_N \in K[B]$
of $M$: $1=\sum_{l=1}^N f_l \Delta_l$ for suitable $f_l \in K[B]$.

Fix a $j \in \{1,\ldots,n\}$. For each $l \in \{1,\ldots,N\}$ we can
write $e_j^T=v_{jl}^T M$, where $e_j$ is the $j$-th standard basis vector
in $K[B]^n$ and $v_{jl}$ is a vector in $K(B)^m$ (supported
only on the positions corresponding to $\Delta_l$) 
satisfying $\Delta_l
v_j \in K[B]^m$; this is just Cramer's rule. Now
\[ e_j=1 e_j=\sum_{l=1}^N f_l \Delta_l e_j=\sum_{l=1}^N f_l (\Delta_l
v_{jl}^T) M, \]
and we conclude that $J_e$ contains an element of the form $k_{0,j} +
x_j^{p^e}$ with $k_{0,j} \in K[B]$.

Consider the morphism $\varphi:A \to A$ dual to the homomorphism $K[A] \to
K[A]$ that restricts to the identity on $K[B]$ and sends each $x_j$ to
$x_j^{p^e}$. Since $K$ is algebraically closed, $\varphi$ is a homeomorphism
in the Zariski-topology, and since $\varphi$ commutes with the projection $A
\to B$, it suffices to show that this projection restricts to a closed
embedding from $Z':=\varphi(Z)$ into $B$. Let $J'$ be the ideal of $Z'$.
By the previous paragraph, $J'$ contains an element of the form
$k_{0,j}
+ x_j$ with $k_{0,j} \in K[B]$ for each $j$. Hence the map $K[B] \to
K[A]/J'$ is surjective, so $Z' \to B$ a closed embedding, as desired.
\end{proof}

\begin{remark}
In characteristic zero, the Zariski-homeomorphism from the lemma is in
fact a closed embedding. In positive characteristic, it need not be.
\end{remark}

\subsection{Extending the field}
Let $P:\VVec \to \VVec$ be a finite-degree polynomial functor over the
infinite field $K$, let $L$ be an extension field of $K$, and denote by
$\VVec_L$ the category of finite-dimensional vector spaces over $L$. We
construct a polynomial functor $P_L:\VVec_L \to \VVec_L$ as follows. For
every $U \in \VVec_L$ we fix a $V_U \in \VVec$ 
and an isomorphism $\psi_U: U \to L \otimes_K V_U$ of
$L$-vector spaces.

At the level of objects, $P_L$ is defined by 
$P_L(U):=L \otimes_K P(V_U)$. To define $P_L$ on morphisms
we proceed as follows. For $U,U' \in \VVec_L$ the
polynomial map $P:\Hom_K(V_U,V_{U'}) \to
\Hom_K(P(V_U),P(V_{U'}))$ extends uniquely to a polynomial map 
\[ P'_L:L \otimes_K \Hom_K(V_U,V_{U'}) \to L \otimes_K
\Hom_K(P(V_U),P(V_{U'})) \]
of the same degree; here we use that $K$ is infinite. The domain and
codomain of $P'_L$ are canonically $\Hom_L(L \otimes_K V_U,L
\otimes_K
V_{U'})$ and $\Hom_L(P_L(U),P_L(U'))$, respectively. Hence for $\varphi
\in \Hom_L(U,U')$ we may set
\[ P_L(\varphi)=P'_L(\psi_{U'} \circ \varphi \circ \psi_U^{-1}). \]
A simple calculation shows that $P_L$ is indeed a
polynomial functor $\VVec_L \to \VVec_L$. 

Furthermore, if $X$ is a $\VVec$-closed subset of $P$, then
we obtain a $\VVec$-closed subset $X_L$ of $P_L$ by letting
$X_L(U)$ be the Zariski closure of $\{1 \otimes q \mid q \in X(V_U)\}$
in $P_L(U)=L \otimes_K P(V_U)$. The
following lemma is straightforward.

\begin{lemma} \label{lm:Field}
The map $X \mapsto X_L$ from $\VVec$-closed subsets of $P$ to
$\VVec_L$-closed subsets of $P_L$ is inclusion-preserving and
injective. Consequently, Noetherianity of $P_L$ implies that
of $P$.
\end{lemma}

\subsection{Proof of Theorem~\ref{thm:Noetherianity}}
\label{ssec:Proof}

If $d=0$, then $P(V)$ is a finite-dimensional space independent of $V$,
and for every linear map $\varphi:V \to W$ the map $P(\varphi)$ is the identity,
so the theorem is the topological corollary to Hilbert's basis theorem. We
therefore may and will assume that $d>0$. Furthermore, by
Lemma~\ref{lm:Field} we may assume that $K$ is algebraically
closed, so that we can use Lemma~\ref{lm:Closed}.

We proceed by induction, assuming that the theorem holds for all
polynomial functors $P' \prec P$ in the well-founded order from
\S\ref{ssec:Order}; this is our {\em outer} induction hypothesis.
From \S\ref{ssec:Splitting} we recall the definition of $\delta_X$.
We now order $\VVec$-varieties inside the fixed $P$ by $X>Y$ if either
$X_Q \supsetneq Y_Q$ or else $X_Q=Y_Q$ and $\delta_X>\delta_Y$.
Since $\delta_X$ takes values in a well-ordered set, for any strictly
decreasing sequence $X_1>X_2>\ldots$, $(X_i)_Q$ must become strictly
smaller infinitely often; but this is impossible since $Q$ is Noetherian
by the outer induction hypothesis.  Hence $>$ is a well-founded order
on $\VVec$-subvarieties of $P$.

We set out to prove, by induction along this well-founded order, that
each $\VVec$-variety $X \subseteq P$ is Noetherian as a $\VVec$-topological
space. Our {\em inner} induction hypothesis states that this holds for
each $\VVec$-variety $Y<X$ inside $P$.

First assume that $\delta_X=\infty$, which means that $X$ is the pre-image
of its projection $X_Q$, and let $Y \subsetneq X$ be any proper closed
$\VVec$-subset. Then either $Y_Q \subsetneq X_Q$ or else $Y_Q=X_Q$ and
$\delta_Y < \delta_X$. Hence $Y < X$, so that $Y$ is Noetherian by the
inner induction hypothesis. Since any inclusion-wise strictly decreasing
chain of closed $\VVec$-subsets of $X$ must contain such a $Y$ as its
first or second element, $X$ is Noetherian, as well.

So we may assume that $\delta_X \in \ZZ_{\geq 1}$.  Take $U \in \VVec$
of minimal dimension for which $\cI_X(U)$ contains a nonzero homogeneous
element of degree $\delta_X$, and let $f \in K[P(U)]$ be a homogeneous
polynomial representing this element. Regarding $f$ as a polynomial with
coefficients from $K[Q(U)]$ in coordinates that restrict to a basis of
$R(U)^*$, we may remove from $f$ all terms with coefficients that vanish
identically on $X_Q(U)$, and then at least one non-constant term
survives. 

By Lemma~\ref{lm:Hasse} applied to $f$ with $W'=P(U)$ and $W=R(U)$,
this implies that there exists an $r_0 \in R(U)$ such that the
directional derivative $h:=\del_{r_0} f \in K[P(U)]$ in the sense of
\S\ref{ssec:Additive} also has at least some coefficient in $K[Q(U)]$
that does not vanish on $X_Q(U)$. Let $e_0 \in \ZZo$ be the exponent
$e$ in the Lemma~\ref{lm:Hasse}; so the map $r \mapsto (\del_r f)(q)$
lies in $\Add(R(U))_{e_0}$ for each $q \in P(U)$. Since coordinate
functions on $R_d(U)$ were assigned degree $d$ (\S\ref{ssec:Var}), we
have $\deg(h)=\deg(f)-d p^{e_0}$ and in particular $\deg(h)<\deg(f)$. By
minimality of the degree of $f$, the polynomial $h$ does not vanish
identically on $X$.

Let $Y$ be the largest closed $\VVec$-subset of $X$ such that $h$ does
vanish identically on $Y(U)$, i.e., $Y(V)=\{p \in X(V) \mid h(P(\varphi)p)=0
\text{ for all } \varphi \in \Hom(V,U)\}$. Then either $Y_Q \subsetneq X_Q$
or else $Y_Q=X_Q$ and $\delta_Y \leq \deg(h) < \delta_X$. Hence $Y<X$,
so $Y$ is Noetherian by the inner induction hypothesis.

Define $Z$ by $Z(V):=X(V) \setminus Y(V)$. This is typically {\em not}
a $\VVec$-subset of $X$, since for $\varphi \in \Hom_\VVec(V,V')$ the map
$P(\varphi)$ might map points of $Z(V)$ into $Y(V')$, i.e., outside $Z(V')$.
Indeed, since we chose $U$ of minimal dimension, when we pull back $f$
to a $P(V)$ for $\dim V<\dim U$, we obtain a polynomial that is zero
modulo the ideal of $X_Q(V)$. This implies that the pull-back of $h$ is
identically zero on $X(V)$, so that $Z(V)=\emptyset$. To remedy this,
we now construct a $\VVec$-variety $Z'$ closely related to
$Z$, by shifting our polynomial functor over $U$ as in
\S\ref{ssec:Shift}.

Set $P':=P \circ \Sh_U$ and $X':=X \circ \Sh_U$ and
consider the open $\VVec$-subset $Z'$ of $X'$ defined by $Z'(V):=\{q \in
X(U \oplus V) \mid h(q) \neq 0\}$.  Here we regard $h$ as a polynomial
on $P(U \oplus V)$ via the map $P(U \oplus V) \to P(U)$ coming from the
projection $\pi_U: U \oplus V \to U$ along $V$. As the maps $\pi_U \circ
g$ with $g \in \GL(U \oplus V)$ are Zariski-dense in $\Hom(U \oplus V,U)$,
$Z'$ and $Z \circ \Sh_U$ are related by
\begin{align} \label{eq:Z} \tag{*}
&Z(U \oplus V)=\{q \in X(U \oplus V) \mid \exists \psi \in
\Hom(U \oplus V,U): h(P(\psi)q) \neq 0\} \\
\notag
&=\{q \in X'(V) \mid h(g(q)) \neq 0 \text{
for some } g \in \GL(U \oplus V) \} = \bigcup_{g \in GL(U \oplus V)} g Z'(V). 
\end{align}
Recall that $R$ is an irreducible subfunctor of $P_d$.  Write $R':=R
\circ \Sh_U=R'_0 \oplus \cdots \oplus R'_d \subseteq P'$ where $R'_e$ is
homogeneous of degree $e$ and $R'_d = R$ by Lemma~\ref{lm:Shift}. Define
$Q':= P' / R'_d$ and note that $Q' \prec P$ since $Q'$ has degree
at most $d$ and the degree-$d$ part of $Q'$ is equal to $P_d/R$
by Lemma~\ref{lm:Shift}. In particular, $Q'$ is a Noetherian
$\VVec$-topological space by the outer induction hypothesis.  

\begin{remark}
Note that for $e<d$, $Q'_e$ typically will have dimension larger than
$P_e$! Also, $Q'$ is not equal to $Q \circ \Sh_U$: there is a surjection
$Q' \to Q \circ \Sh_U$ with kernel $R'/R'_d$.
\end{remark}

The map $P(\pi_U): P(U \oplus V) \to P(U)$ has $R'_d(V)$ in its kernel, so we can regard $h$
as a polynomial on $Q'(V)$, as well. Consider the basic open
$\VVec$-subset
of $Q'$ defined by $B(V):=\{q \in Q'(V) \mid h(q) \neq 0\}.$ As $Q'$
is Noetherian, so is $B$ with its induced topology.

\begin{lemma} \label{lm:Closed2}
For every $K$-vector space $V$ the projection $Z'(V) \to B(V)$ is a
Zariski-homeomorphism with a closed subset of $B(V)$.
\end{lemma}

Before proving this lemma, we use it to complete the proof of
Theorem~\ref{thm:Noetherianity}. First, since $B$ is Noetherian,
Lemma~\ref{lm:Closed2} implies that so is $Z'$.

Then suppose that $X=X_0 \supseteq X_1 \supseteq \cdots$ is a sequence
of closed $\VVec$-subsets of $X$. By Noetherianity of $Y$ there exists an
$n_0$ such that for all $V \in \VVec$ the sequence $(X_n(V) \cap Y(V))_n$
is constant for $n \geq n_0$. By Noetherianity of $Z'$ there exists
an $n_1$ such that for all $V \in \VVec$ the sequence $(X_n(U \oplus V)
\cap Z'(V))_n$ is constant for $n \geq n_1$. Using \eqref{eq:Z} we find
\begin{align*} X_n(U \oplus V) \cap Z(U \oplus V)
&=
\bigcup_{g \in \GL(U \oplus V)} X_n(U \oplus V) \cap gZ'(V)\\
&= 
\bigcup_{g \in \GL(U \oplus V)} g(X_n(U \oplus V) \cap Z'(V)), 
\end{align*}
where in the last step we used that $X_n(U \oplus V)$ is $\GL(U \oplus
V)$-stable. We find that, for each $V$ of dimension at least $\dim U$, 
the sequence $(X_n(V) \cap Z(V))_n$ is constant for $n \geq n_1$.
Since $Z(V)=\emptyset$ for $V$ of dimension less than $\dim U$, we
find that $X_n$ is constant for $n \geq \max\{n_0,n_1\}$. This proves
Noetherianity of $X$ and concludes the proof of the inner induction step.
$\hfill \square$

\begin{example} \label{ex:Running4}
We pause and see what Lemma~\ref{lm:Closed2} says in the running
Example~\ref{ex:Running},\ref{ex:Running2},\ref{ex:Running3}.
Here $U=K^2$ and we take $V=K^n$. Then 
$Z'(V)$ is the variety of $(n+2) \times (n+2)$-matrices
\[ M=\begin{bmatrix} C & D \\ E & F
\end{bmatrix} \]
of rank $1$ such that the upper-left $2 \times 2$-submatrix $C$ has a
non-zero coordinate $z_{12}$ (recall that $h=2z_{12}$). 

The lemma says
that forgetting the skew-symmetric part of $F$ is a closed embedding
of $Z'(V)$ into the open subset of $U \otimes U + V \otimes U + U \otimes
V + Q(V)$ where $h$ is nonzero. We prove this by showing that, on $Z'$,
each entry of the skew-symmetric part of $F$ can be expressed as a
rational function in the entries of $C,D,E$ and the entries
of the symmetric part of $F$, with a denominator equal to $h$.

We assume $n \geq 2$ and consider a surjective linear
map $\varphi:V \to U$. Thinking of $\varphi$ as a $2 \times n$-matrix, we have
\[ P(1_U \oplus t\varphi)M= C + t(\varphi E + D \varphi^T) + t^2 \varphi F \varphi^T. \]
Since $Z'$ is a $\VVec$-closed subset, for all $M \in Z'(V)$, the
$2 \times 2$-determinant $f$ vanishes on the latter $2 \times 2$-matrix
for all choices of $t$ and of $\varphi$. Hence expanding $f(P(1_U
\oplus t \varphi)M)$
as a polynomial in $t$, the coefficient $k$ of $t^2$ also vanishes for
all $M,\varphi$. Take $i<j$ in $\{1,\ldots,n\}$ and specialise $\varphi$ to the linear
map sending the $i$-th standard basis vector to $e_1$, the $j$-th standard
basis vector to $e_2$, and all other standard basis vectors
to zero. Then the matrix above reads
\[ \begin{bmatrix} 
C_{11} + t E_{i1} + t D_{1i} + t^2 F_{ii} & 
C_{12} + t E_{i2} + t D_{1j} + t^2 F_{ij} \\
C_{21} + t E_{j1} + t D_{2i} + t^2 F_{ji} &
C_{22} + t E_{j2} + t D_{2j} + t^2 F_{jj}
\end{bmatrix}.
\]
The coefficient of $t^2$ in the determinant of
this matrix is
\[ k(C,D,E,F)=F_{ii} C_{22} + F_{jj} C_{11} - 
F_{ij} C_{21} - F_{ji} C_{12} + \cdots \]
where the remaining terms do not involve $F$. Using, as in
Example~\ref{ex:Running}, the variables $y,z$ for the
symmetric and skew-symmetric parts of $C$, and using the variables
$y',z'$ for the symmetric and skew-symmetric parts of $F$,
this reads
\begin{align*}  &y'_{ii} y_{22} + y'_{jj} y_{11} 
-(y'_{ij}+z'_{ij})(y_{12}-z_{12}) 
-(y'_{ij}-z'_{ij})(y_{12}+z_{12})
+ \cdots\\ &= 2 z'_{ij} z_{12} + \cdots = z'_{ij} h(C) +
\cdots  \end{align*}
where the dots in the last two expressions contain only variables that we do
not discard in the projection $P'(V) \to Q'(V)$. This shows
that, on $Z'$, the coordinate $z'_{ij}$ can be expressed in
the entries of $C,D,E$ and the coordinates $y'_{ij}$ on the
symmetric part of $F$, as desired. \hfill $\clubsuit$
\end{example}

\begin{proof}[Proof of Lemma~\ref{lm:Closed2}.]
First assume that $\dim V \geq \max\{\dim U,d\}$. We want to apply
Lemma~\ref{lm:Closed} with $W'$ equal to $P'(V)$, $W$ equal to
$R'_d(V)=R(V)$, $B$ equal to $B(V)$, and $Z$ equal to $Z'(V)$. Let
$J$ be the ideal of $Z'(V)$ in the coordinate ring of the pre-image
of $B(V)$ inside $P'(V)$, and let $J_{e_0}$ be as in the text preceding
Lemma~\ref{lm:Closed}.

Fix any surjective linear map $\varphi:V \to U$.
For $t \in K$ and $e \in \{0,\ldots,d\}$ consider the linear map
$\Phi_e(t):=P_e(1_U \oplus t\varphi): P_e(U \oplus V) \to P_e(U)$. This map
depends on $t$ as a polynomial of degree at most $e$, hence decomposes
as $\Phi_e(t)=t^0 \Phi_{e0} + \cdots + t^e \Phi_{ee}$ for unique linear
maps $\Phi_{ei}:P_e(U \oplus V) \to P_e(U)$. Note that
$\Phi_{ee}=P_e(0 \oplus \varphi)$ and $\Phi_{e0}=P_e(1_U
\oplus 0)=P_e(\pi_U)$.  

On the other hand,
decompose $P'_e:=P_e \circ \Sh_u=P'_{e0} \oplus \cdots \oplus P'_{ee}$
where $P'_{ei}$ is a homogeneous polynomial functor of degree $i$. Then $\Phi_{ei}$ is zero on
$P'_{ej}(V)$ except when $i=j$. To see this, take a $q \in
P'_{ej}(V)$, compute 
\begin{align*} \Phi_e(t)q&=P_e(1_U \oplus t\varphi)q= P_e(1_U \oplus \varphi)
P_e(1_U \oplus t1_V) q\\&= P_e(1_U \oplus \varphi) P'_e(t 1_V)q= t^j P_e(1_U
\oplus \varphi) q, \end{align*}
and observe that the right-hand side is homogeneous
of degree $j$ in $t$.

The $\Phi_e(t)$ together form the map $\Phi(t):=P(1_U \oplus t\varphi): P(U
\oplus V) \to P(U)$. Since $X$ is a $\VVec$-subset of $P$
and $f$ vanishes on $X(U)$, $f(\Phi(t)q)=0$ for all $q \in
X(U \oplus V)$ and $t \in K$. 
This implies that the
coefficient $k(q)$
of $t^{dp^{e_0}}$ in $f(\Phi(t)q)$ also vanishes identically on $X(U
\oplus V)$. 

To determine how $k(q)$ depends on $R_d'(V)$,
consider $r \in R'_d(V)$ and $q=q_0 + \cdots + q_d \in P(U \oplus V)$
with $q_e \in P_e(U \oplus V)$, and for variables $t,s$ compute
\begin{align*} 
&f\left(\Phi(t) (q_0 + \cdots + q_{d-1} + (q_d+s r))\right)\\
&= f(\Phi_0(t) q_0 + \cdots + 
\Phi_{d-1}(t) q_{d-1} + \Phi_d(t)(q_d+s r))\\
&= f(\Phi_0(t) q_0 + \cdots + \Phi_{d-1}(t) q_{d-1}
+\Phi_d(t)q_d+t^d s \Phi_{dd} r)\\
&= f(\Phi_0(t) q_0 + \cdots + \Phi_{d-1}(t) q_{d-1} 
+\Phi_d(t)q_d+t^d s P_d(0 \oplus \varphi) r)\\
&\equiv f(\Phi_0(t) q_0 + \cdots + \Phi_{d-1}(t) q_{d-1} 
+\Phi_d(t)q_d)\\
&+(t^d s)^{p^{e_0}} (\partial_{P_d(0
\oplus \varphi)r} f)(\Phi_0(t) q_0 + \cdots + \Phi_{d-1}(t)
q_{d-1}
+\Phi_d(t)q_d) \mod ((t^d s)^{p^{e_0}+1}) 
\end{align*}
where in the second equality we have used that $r \in P_{dd}'(V)$
and in the last step we have used Lemma~\ref{lm:Hasse}. We see that
$t^{dp^{e_0}}$ is the lowest power of $t$ to whose coefficient $r$
contributes, and that this contribution is additive in $r$. More
specifically, for each $s \in K$ we have $k(q+s r)=k(q)+s^{p^{e_0}}
(\partial_r k)(q)$, and $(\partial_r k)(q)$ equals the value of the
directional derivative $\del_{R (\varphi)r} f$ at the point $\sum_{e=0}^d
\Phi_{e0} q_e=P(\pi_U)(q)$. In particular, $k \in J_{e_0}$.

From now on, assume that the image of $q$ in $Q'(V)$ lies in $B(V)$. Since
$\varphi$ is surjective, so is $R(\varphi):R(V) \to R(U)$. In particular, there
exists an $r \in R(V)$ such that $R(\varphi)r=r_0$.  For such an $r$ we have
$\ell(r):=(\partial_r k)(q)=h(q) \neq 0$, so $\ell \in \Add(R(V))_{e_0}$
is not zero.

Keeping $q$ fixed but replacing $\varphi$ by $\varphi \circ g$ for
$g \in \GL(V)$, $\ell$ transforms into the additive function $r \mapsto
\ell(R(g)r)$. Hence by varying $g$ we find that the image of $J_e$ in
$\Add(R(V))_e$ under the map $\tilde{k} \mapsto (r \mapsto (\partial_r
\tilde{k})(q))$ from Lemma~\ref{lm:Closed} contains a nonzero a nonzero
$\GL(V)$-submodule $L$ of $\Add(R(V))_e$. Since $R$ is irreducible and
$\dim V \geq d$, $R(V)$ is an irreducible $\GL(V)$-module \cite[Lemma
3.4]{Friedlander97}, and this implies the irreducibility of $(R(V))^*$
and of $\Add(R(V))_{e_0}$---indeed, raising to the power $p^{e_0}$ gives
a bijection from $\GL(V)$-submodules to $\GL(V)$-submodules of the latter.

We conclude that $L=\Add(R(V))_{e_0}$, and since $q$ was arbitrary in
the pre-image of $B(V)$, the conditions of Lemma~\ref{lm:Closed}
are fulfilled. Hence the projection $Z'(V) \to B(V)$ is a
Zariski-homeomorphism with a closed subset of $B(V)$, as desired.

Finally, if $\dim V < \max \{\dim U,d\}$, then take any embedding $\iota:
V \to V'$ where $V'$ does have sufficiently high dimension.
Then we have a commuting diagram 
\[ 
\xymatrix@R-1pc{
Z'(V) \ar[d] \ar[r]^{P'(\iota)} & Z'(V') \ar[d] \\
B(V) \ar[r]_{P'(\iota)} & B(V')
}
\]
where all arrows except, {\em a priori}, the left-most one are
homeomorphisms with closed subsets of the target space. But then so is
the left-most one.
\end{proof}

\subsection{Comments on the proof}
\label{ssec:Comments}

The idea to do induction on $P$ is not completely new: it is also
used, in special cases, in \cite{Draisma14a, Eggermont14, Derksen17,
Hochster16}.  In these papers, more information than just Noetherianity
is extracted from the proof: e.g.~that the {\em tuple rank} of a tuple
of matrices is bounded in a proper closed subvariety of a polynomial
functor capturing matrices \cite{Draisma14a,Eggermont14}, or that the
{\em $q$-rank} of a cubic is bounded \cite{Derksen17}, or that the {\em
strength} of a homogeneous form is bounded \cite{Hochster16}. Our proof
above does not directly yield such qualitative information. However,
in \cite{Draisma18b} we repair this defect for symmetric, alternating,
and ordinary tensors and characteristic zero or sufficiently large.

If $K$ has characteristic zero, then symmetric powers $S^d V$ are
irreducible $\GL(V)$-modules, and one can prove Noetherianity for
direct sums of these without the need for more general polynomial
functors---though also without the proof becoming any easier.
But in general characteristic, symmetric powers need not be
irreducible, and polynomial functors need not be completely reducible
into irreducible summands, so reducing modulo an irreducible subfunctor
is the only natural thing to do.

The idea further to do induction on $\delta_X$ and to use directional
derivatives {\em is} new, but inspired by techniques used earlier in
\cite{Draisma08b,Draisma11d,Draisma14a}, where a determinant is regarded
as an affine-linear polynomial in one matrix entry, whose coefficient
is a determinant of lower order, and induction is done over
that order.

\subsection{An open problem}

The most tantalising open problem in this area is the following.

\begin{question}
Let $P$ be a finite-degree polynomial functor over an infinite field
$K$. Does any sequence $I_1 \subseteq I_2 \subseteq \cdots$ of ideals
in $K[P]$ eventually become constant?
\end{question}

For $P$ of degree at most $1$, the answer is yes, and it follows from
the stronger statement that the ring $R[y_{ij} | i=1,\ldots,k,
j \in \NN]$, acted upon by $\Sym(\NN)$ via $\pi y_{ij}=y_{i \pi(j)}$,
is $\Sym(\NN)$-Noetherian for any Noetherian ground ring $R$
\cite{Aschenbrenner07,Cohen87,Hillar08}. For $P=S^2$ and $P=\Wedge^2$
in characteristic zero, the anwer is also yes, since we know all
$\GL(V)$-stable ideals from \cite{Abeasis80a} and \cite{Abeasis80b},
respectively. In \cite{Nagpal15} a much stronger result than this is
established for $S^2$ and $\Wedge^2$, namely, that finitely generated modules
over $K[P]$ are also Noetherian. These questions were first raised,
in the setting of twisted commutative algebras, in \cite{Snowden10}.
They remain widely open even for more general degree-two functors,
and also for, say, $\Wedge^2$ in positive characteristic.

\ \\
\ \\
\bibliographystyle{amsplain}
\bibliography{draismajournal,draismapreprint,diffeq}

\end{document}